\documentclass[12pt]{article}
\usepackage{amsmath}
\usepackage{amsfonts}
\usepackage{amsthm}
\usepackage{times}
\def\acts{\triangleright}
\def\Z{{\mathbb Z}}
\def\C{{\mathbb C}}
\def\R{{\mathbb R}}

\def\CA{{\mathcal A}}

\def\CL{{\mathcal L}}

\def\CV{{\mathcal V}}

\DeclareMathSymbol\crossrt{\mathrel}{AMSb}{"6E}
\DeclareMathSymbol\crosslt{\mathrel}{AMSb}{"6F}

\def\id{{\mathrm i}{\mathrm d}}
\def\ts{\otimes}

\def\oh{\frac{1}{2}}

\def\acts{\triangleright}

\def\ts{\otimes}
\def\ket#1{| #1 \rangle}

\newtheorem{lemma}{Lemma}[section]
\newtheorem{theorem}[lemma]{Theorem}

\def\half{{\frac{1}{2}}}

\title{On Spin Structures and Dirac Operators \\ on the Noncommutative Torus}
\author{ Mario Paschke \\
{\em Max-Planck-Institut f\"ur Mathematik in den Naturwissenschaften} \\
{\em Inselstr.\ 22, 04103 Leipzig, Germany} \\
\and Andrzej Sitarz\thanks{Alexander von Humboldt Fellow}
\thanks{Partially supported by MNII Grant 115/E-343/SPB/6.PR UE/DIE 50/2005--2008}\\
{\em Institute of Physics, Jagiellonian University}\\
{\em Reymonta 4, 30-059 Krak\'ow, Poland} \\
and \\
{\em Mathematisches Institut, Heinrich-Heine-Universit\"at}\\
{\em Universit\"atsstrasse 1, 40225 D\"usseldorf, Germany}
}
\begin{document}
\maketitle
\begin{abstract}
\noindent We find and classify possible equivariant spin
structures with Dirac operators on the noncommutative torus,
proving that similarly as in the classical case the spectrum of
the Dirac operator depends on the spin structure.
\end{abstract}

\section{Introduction}

Unlike classical differential geometry, in Noncommutative Geometry
there seems to be no well-established, commonly accepted notion of
spin structures over noncommutative
manifolds. The notion of a real spectral triple \cite{Connes95}, which
provides a noncommutative counterpart of the spin manifold definition
incorporates in the classical case both spin structure and the Dirac
operator. In the known examples, it is a priori difficult to separate
these two ingredients. On the other hand, one would expect that possible
spectral triple constructions should - like in the classical case - depend
on some possible choices of the reality structure $J$. The resulting
Dirac operators might then have different spectra.

In the paper we investigate the simplest and most studied example
of the two-dimensional noncommutative torus. In the classical limit,
it is well known that there are four different spin structure on the
torus and the Dirac operators have different spectra \cite{Baer}. We
recover this result in the noncommutative case, showing that out of
several {\em a priori} admissible reality structures only some admit
reasonable Dirac operators.

For the details of the spectral triple axiomatics we refer the reader
to \cite{GBVF}. In our approach we use equivariance under the
global symmetry of the noncommutative manifolds as defined in
\cite{MP,Sitarz03}.

\section{Spectral triples on Noncommutative Torus}

\subsection{Equivariant representation and a real structure}

Let $U,V$ be  unitary elements generating the algebra of
polynomial functions on the noncommutative torus $\CA(T_\lambda)$,
$\lambda$ a generic complex number, $|\lambda|=1$,
\begin{equation}
UV = \lambda VU,
\label{nctorus}
\end{equation}
and $\delta_1, \delta_2$ be the basis of derivations acting
on $\CA(T_\lambda)$, which represent the Lie algebra type symmetry
(denoted later by $\CL$) of the noncommutative torus:
\begin{equation}
\begin{split}
\delta_1 \acts U &= U, \quad\quad\quad \delta_2 \acts U = 0, \\
\delta_1 \acts V &= 0, \quad\quad\quad \delta_2 \acts V = V.
\end{split}
\end{equation}

We look for all possible equivariant representations with
an equivariant real structure $J$ and, since the spectral triple for
a two-dimensional torus should be even, for a $\Z_2$ grading $\gamma$.
Note that in this case the Hilbert space must be graded and we need
to find two equivariant representations $\pi_1$, $\pi_2$.

First, we assume that $\CV$ is a vector space on which we have
a well-defined star-representation $\rho$ of the derivations
and that $e_{\mu,\nu}$ are their mutual eigenvectors:
\begin{align}
\rho(\delta_1) e_{\mu,\nu} &= \mu \, e_{\mu,\nu}, &
\rho(\delta_2) e_{\mu,\nu} &= \nu \, e_{\mu,\nu},
\end{align}
where $\mu,\nu \in \R$ are arbitrary numerical labels.

The equivariance condition, for $l \in \CL$ and
$a \in \CA(T_\lambda)$, $v \in \CV$:
\begin{equation}
 \rho(\ell) \pi(a) v = \left( \pi(\ell \acts a) + \pi(a) \right) v,
\end{equation}
gives us the following result
\begin{align}
\pi(U) e_{\mu,\nu} &= u_{\mu,\nu} \, e_{\mu+1,\nu},&
\pi(V) e_{\mu,\nu} &= v_{\mu,\nu} \, e_{\mu,\nu+1}.
\end{align}

Further, from the commutation relation (\ref{nctorus}) we get
$$ u_{\mu,\nu+1} v_{\mu,\nu} =
\lambda \, u_{\mu,\nu} v_{\mu+1,\nu}. $$
Since the eigenvectors $e_{\mu,\nu}$ can always be rescaled by
a suitable phase we might choose,
\begin{equation}
u_{\mu,\nu} = 1, \;\;\;\;\;\; v_{\mu,\nu} = \lambda^{-\mu}.
\end{equation}

Therefore, the minimal irreducible equivariant representation
would consists of a linear span of $e_{\mu+m,\nu+n}$ for all $m,n \in \Z$:
\begin{equation}
\CV_{\mu,\nu} = \bigoplus_{m,n \in \Z} V_{\mu_0+m,\nu_0+n},
\end{equation}
where each $V_{\mu_0+m,\nu_0+n} \simeq \C$.

Note that we have obtained no restriction for the values of $\mu_0,\nu_0$
and, {\em a priori}, the values for each of the two representations
might be chosen independently and be completely different. We shall
label the two eigenspaces of the grading $\CV^+$ and $\CV^-$ and take
$\CV = \CV_+ \oplus \CV_-$.

Next, we look for an antilinear unitary operator $J$ on $\CV$, such that:
\begin{equation}
J^2=-1, \;\;\;\;\; J\gamma = - \gamma J.
\label{Jcond}
\end{equation}
The latter condition means that $J$ maps $\CV_\pm$ onto $\CV_\mp$.
We use the equivariance condition for $J$, which, in the case of
Lie-algebra symmetry, reads:
\begin{equation}
- J \rho(\ell^*) v = \rho(\ell) J v, \;\; \forall v \in \CV, \ell \in \CL.
\end{equation}
Writing it explicitly for the derivations $\delta_1,\delta_2$:
\begin{equation}
\begin{aligned}
 - (\mu_1+m_1) J \, e_{\mu_1+m_1,\nu_1+n_1, \pm} &=
\delta_1 \left( J e_{\mu_1+m_1,\nu_1+n_1, \pm} \right), \\
- (\nu_1+n_1) J e_{\mu_1+m_1,\nu_1+n_1, \pm} &=
\delta_2 \left( J e_{\mu_1+m_1,\nu_1+n_1, \pm} \right).
\label{Jeq0}
\end{aligned}
\end{equation}
Hence, we immediately get:
\begin{lemma}
An equivariant $J: \CV \to \CV$ satisfying (\ref{Jcond}) exists only
and only if the spectrum of the derivations $\delta_1,\delta_2$ on
spaces $\CV_+$ is a symmetric image (with respect to $x \to -x$) of
the spectrum on $\CV_-$. Then:
\begin{equation}
J e_{\mu,\nu,\pm} = \pm j(\mu,\nu) \, e_{-\mu, -\nu, \mp},
\label{Jeq}
\end{equation}
where $j(\mu,\nu)$ is an arbitrary phase and $e_{\mu,\nu,\pm} \in \CV_\pm$.
\end{lemma}
\begin{proof}
The equation (\ref{Jeq}) is direct consequence of (\ref{Jeq0}).
Clearly, if $\mu$ is in the spectrum of $\delta_1$ on $\CV_+$
then $-\mu$ must be in the spectrum of $\delta_1$ on $\CV_-$.
\end{proof}
Next, by requiring that $J$ maps the algebra to its commutant, we have
\begin{lemma}
$J$ maps the algebra $\CA(T^2_\lambda)$ to its commutant
if and only if
$$ j(\mu,\nu) =e^{i \phi \mu + i \psi \nu + \theta} \lambda^{-\mu\nu}, $$
where $\phi,\psi$ and $\theta$ are arbitrary real numbers.
\end{lemma}
\begin{proof}
First, observe that $J$ is of the form $J= \Lambda J_0$, where
$J_0$ is antilinear and $\Lambda$ is diagonal and unitary:
$$  J_0 e_{\mu,\nu,\pm} = \pm e_{-\mu,-\nu,\mp}, \;\;\;\;\;\;\;
\Lambda e_{\mu,\nu,\pm} = j(-\mu,-\nu) e_{\mu,\nu,\pm}.$$

Taking first $j(\mu,\nu) = \lambda^{-\mu\nu}$ we verify
explicitly that such canonical $J_c$ maps the algebra to the commutant:
\begin{align}
U^o e_{\mu,\nu,\pm} = J_c^* U^* J_c e_{\mu,\nu,\pm} &=
 \lambda^{-\nu} e_{\mu+1,\nu,\pm}, \\
V^o e_{\mu,\nu,\pm} = J_c^* V^* J_c e_{\mu,\nu,\pm} &=
 e_{\mu,\nu+1,\pm},
\end{align}
and it is easy to verify that $U^o,V^o$ are indeed in
the commutant of $\CA(T_\lambda)$.

Now, assume that there exists a different diagonal unitary
$\Lambda'$, such that $J' = \Lambda J_0$ maps the algebra
to the commutant. Consider the map $W = J' J_c$:
$$ J' J_c e_{\mu,\nu,\pm} = w_\pm(\mu,\nu) e_{\mu,\nu,\pm}.$$
Clearly, $W = J' J_c$ is unitary. Both $W^* U W$, and $W^* V W$
should commute with $U^o$ and $W^o$, for instance:
$$ \left[(J' J_c)^* U (J' J_c), J_c^* U J_c \right]
= J_c^* \left[ (J')^* U J', U \right] J_c  = 0, $$ since we assumed
that $J'$ also maps $\CA(T_\lambda)$ to the commutant.
Let us calculate $[ W^* U W, V^o]$:
\begin{equation}
\begin{split}
\left[W^* U W, U^o \right] e_{\mu,\nu,\pm} =&
 \lambda^{-\nu} \left( w_\pm(\mu+1,\nu) w_\pm(\mu+2,\nu)^* \right. \\
 &- \left. w_\pm(\mu,\nu) w_\pm(\mu+1,\nu)^* \right) e_{\mu+2,\nu,\pm}.
\end{split}
\end{equation}
Therefore, using the unitarity of $W$:
$$ \left(w_\pm(\mu+1,\nu)\right)^2 =  w_\pm(\mu,\nu) w_\pm(\mu+2,\mu), $$
This recurrence relation has the following solution:
\begin{equation}
w(\mu,\nu) = e^{i \mu \phi_\pm} w_\pm^0(\nu),
\label{recur-1}
\end{equation}
where $\phi_\pm$ are arbitrary constants and $w_\pm^0$ an arbitrary
unitary function of $\nu$. Next:
\begin{equation}
\begin{split}
\left[W^* V W, V^o \right] e_{\mu,\nu,\pm} =&
= \lambda^{-\mu} \left( w_\pm(\mu,\nu) w_\pm(\mu,\nu+1)^* \right. \\
 &- \left. w_\pm(\mu,\nu+1) w_\pm(\mu,\nu+2)^* \right) e_{\mu,\nu+2,\pm}.
\end{split}
\end{equation}
Using (\ref{recur-1}) this leads to:
\begin{equation}
\left(w_\pm^0(\nu+1)\right)^2 = w_\pm^0(\nu) w_\pm^0(\nu+2),
\end{equation}
which yields the final general solution:
\begin{equation}
w_\pm(\mu,\nu) = e^{i \phi_\pm \mu + i \psi_\pm \nu + i \theta_\pm}.
\end{equation}
It is an easy exercise to see that the remaining commutators
vanish as well.

The relation between $\phi_+$ and $\phi_-$ (and similarly for
$\psi_\pm,\theta_\pm$) can be fixed using the requirement of
unitarity of $J'$ and the first relation (\ref{Jcond}). Hence,
one gets that: $J' = - W J$ must be:
\begin{equation}
J' e_{\mu,\nu,\pm} =
\pm e^{\pm (i\phi\mu + i \psi\nu + i \theta)}
\lambda^{-\mu\nu} e_{-\mu,-\nu,\mp}.
\end{equation}
The constant $\theta$ is a global phase and can be fixed to $\theta=0$.
\end{proof}

It is interesting to observe what the map $W$ is doing. Clearly
$W^* \pi(a) W$ commutes with the commutant of $\CA(T_\lambda)$ but
is not in the chosen representation of $\CA(T_\lambda)$. In fact,
the conjugation by $W$ induces an automorphism of $\CA(T_\lambda)$
but only when restricted to the single representation $\pi_+$
or $\pi_-$:
\begin{equation}
\begin{aligned}
W^*_+ \pi_+(U) W_+ &= \pi_+( e^{-i \phi} U ), &
W^*_- \pi_-(U) W_- &= \pi_-( e^{i \phi} U ), \\
W^*_+ \pi_+(V) W_+ &= \pi_+( e^{-i \psi} V ), &
W^*_- \pi_-(V) W_- &= \pi_-( e^{i \psi} V ),
\end{aligned}
\end{equation}
where $W_\pm$ denotes the respective diagonal components of $W$.

\subsection{The equivariant Dirac operator}

We proceed now with the construction of an equivariant Dirac operator.
Assume the existence of an equivariant graded linear operator, that is
an operator, which commutes with derivations and anticommutes with
$\gamma$. {}From this we infer that:
\begin{equation}
D e_{\mu,\nu,\pm} = d_{\mu,\nu,\pm} \, e_{\mu,\nu,\mp}.
\end{equation}

Therefore, $D$ intertwines  vectors of the
same eigenvalues of the derivations $\delta_1,\delta_2$  and
in both subspaces such vectors must be present. On the other
hand, from the action of $J$ we know that for each vector
with eigenvalues $\mu,\nu$ in one space there exists one with
eigenvalues $-\mu,-\nu$. Hence, we must have:
\begin{equation}
\exists m,n \in \Z: \mu=-\mu+m, \; \nu=-\nu+n.
\end{equation}
so $\mu, \nu$ are either integers or half-integers.

We summarize the result
\begin{lemma}
There are four possible classes of inequivalent real spectral
triples over the noncommutative torus, given by the following
data. Let $\epsilon_1,\epsilon_2 \in \{0, \half\}$.
Then $\CV_+ = \CV_-$ is a linear span of the orthonormal
vectors $e_{\mu,\nu}$ labelled by $(\mu,\nu) \in
(\Z+\epsilon_1, \Z+\epsilon_2)$, the grading $\gamma$ being:
\begin{equation}
\gamma e_{\mu, \nu, \pm} = \pm e_{\mu, \nu, \pm},
\end{equation}
the real structure $J$:
\begin{equation}
J e_{\mu, \nu, \pm} = e^{\pm i(\phi\mu+\psi\nu+\theta)}
\lambda^{-\mu\nu} e_{-\mu, -\nu, \mp},
\end{equation}
and the equivariant hermitian Dirac operator $D$,
\begin{equation}
D e_{\mu, \nu, \pm} = d_{\mu,\nu}^\pm e_{\mu, \nu, \mp},
\end{equation}
which is determined by the order-one condition.
\end{lemma}

\begin{proof}
First let us prove the existence and uniqueness of an equivariant
$D$ satisfying the order-one condition. From the fact that
$D=D^\dagger$ we learn
$$d_{\mu,\nu}^+ = (d_{\mu,\nu}^-)^*.$$
The order one condition gives the equations:
\begin{align}
(d_{\mu+1,\nu}^+ - d_{\mu,\nu}^+) e^{i\phi} &=
   e^{-i\phi} (d_{\mu,\nu}^+ -d_{\mu-1,\nu}^+), \\
(d_{\mu+1,\nu}^+ - d_{\mu,\nu}^+) e^{i\psi} &=
   e^{-i\psi} (d_{\mu+1,\nu-1}^+ - d_{\mu,\nu-1}^+), \\
(d_{\mu,\nu+1}^+ - d_{\mu,\nu}^+) e^{i\phi} &=
   e^{-i\phi} (d_{\mu-1,\nu+1}^+ - d_{\mu-1,\nu}^+), \\
(d_{\mu,\nu+1}^+ - d_{\mu,\nu}^+) e^{i\psi} &=
   e^{-i\psi} (d_{\mu,\nu}^+ - d_{\mu,\nu-1}^+),
\end{align}
whose only solutions are:
\begin{equation}
d_{\mu,\nu}^+ = \left\{
\begin{aligned}
\tau_0 e^{-2i \phi\mu - 2i \psi\nu} + \epsilon, & \qquad\hbox{when\ } \phi\not= 0, \psi\not=0, \\
\tau_1 \mu + \tau_0 e^{- 2i \psi\nu} + \epsilon, &\qquad \hbox{when\ } \phi= 0, \psi\not=0, \\
\tau_2 \nu + \tau_0 e^{-2i \phi\mu} + \epsilon, & \qquad\hbox{when\ } \phi\not= 0, \psi=0, \\
\tau_1 \mu + \tau_2 \nu + \epsilon. & \qquad\hbox{when\ } \phi= 0, \psi=0,
\end{aligned} \right.
\end{equation}

Finally, demanding that $JD=DJ$ we have, in each of the possible cases:
\begin{itemize}
\item $\phi\not= 0, \psi\not=0$
$$ -(\tau_0)^* e^{i\phi\mu+i\psi\nu} - \epsilon^* e^{-i\phi\mu-i\psi\nu}
 =   (\tau_0)^* e^{-i\phi\mu-i\psi\nu} + (\epsilon)^* e^{i\phi\mu+i\psi\nu},$$
 which holds if $\tau_0 = - \epsilon$,
\item $\phi= 0, \psi\not=0$
$$ -(\tau_1)^*\mu  - (\tau_0)^* e^{i\psi\mu} - \epsilon^* e^{-i\psi\mu}
 =  -(\tau_1)^*\mu + \tau_0 e^{-i\psi\mu} + \epsilon e^{i\psi\mu}, $$
which holds for arbitrary $\tau_1$ and $\tau_0 = - \epsilon$,
\item $\phi\not= 0, \psi=0$
$$ -(\tau_2)^*\nu -(\tau_0)^* e^{i\phi\mu} - \epsilon^* e^{-i\phi\mu}
 =   -(\tau_2)^*\nu + \tau_0 e^{-i\phi\mu} + \epsilon e^{i\phi\mu}, $$
which holds for arbitrary $\tau_2$ and $\tau_0 = - \epsilon$,
\item $\phi\not= 0, \psi\not=0$
$$ -(\tau_1)^*\mu  - (\tau_2)^* \nu  - \epsilon^*
 =  -(\tau_1)^*\mu + - (\tau_2)^* \nu  + \epsilon^*,$$
which holds for arbitrary $\tau_1,\tau_2$ and $\epsilon=0$.
\end{itemize}
\end{proof}

The necessary final condition to fix the spectral data comes either
{}from the requirement of the spectral properties of the Dirac operator
or the Hochschild cycle condition. In the first case, it is obvious
that only for $\phi=\psi=0$, the growth of eigenvalues of the Dirac
operator corresponds to the required axioms and $D$ has compact
resolvent. The Hochschild cocycle condition is more complicated,
and we prove it next.
\begin{lemma}
The Hochschild cycle condition, that is, that there exists:
$$ c = \sum c_0 \ts c_0^o \ts c_1 \ts c_2 \in Z_2(\CA, \CA \ts \CA^o),$$
such that
$$ \gamma = \pi(c) = \pi(c_0) (J^{-1} \pi(c_0^o) J) [ D, \pi(c_1)] [D,\pi(c_1)],$$
is satisfied only if $\phi=\psi=0$.
\end{lemma}
\begin{proof}
To prove it, we start with the case $\phi \not= 0$, $\psi \not= 0$
and observe,
$$ [D,\pi(U)]_\pm = (e^{\mp 2i \phi} -1) U D'_\pm, \;\;\;
[D,\pi(V)]_\pm = (e^{\mp 2i \psi} -1) U D'_\pm, $$
where the sign denotes the restriction of the operators to $\CV_\pm$,
and $D' = D + \tau_0$. Moreover:
$$ D' \pi(U) = e^{\mp 2i \phi} \pi(U) D', \;\;\;
D' \pi(V) = e^{\mp 2i \psi} \pi(v) D'. $$
Any expression of the type:
$ \pi(a_0) [D, \pi(a_1)] [D, \pi(a_2)]$, where $a_0,a_1,a_2$
are homogeneous polynomials in $U,V$ (since we are working with
the algebra of polynomials, we can always restrict ourselves to this
case), must therefore be proportional to:
$$ C(a_1,a_2) \pi(a_0) \pi(a_1) \pi(a_2), $$
where $C(a_1,a_2)$ is a complex number depending only on multi-degree
of polynomials $a_1,a_2$ and the $\psi$ and $\phi$.

Now, assuming that the cocycle condition holds, we would have
a decomposition of $1$ (when restricted to $\CV_+$) as a finite sum
of homogeneous polynomials in $U,V, U^o, V^o$. This, however, is not
possible, unless the polynomials are all of degree $0$. Hence $c_0^o$
{}from the cocycle must be $1$.

Assume next that the cocyle $c$ is trivial, i.e. $c = bc'$. It can
be easily verified that in such case its image $\pi(c)$ is a sum
of commutators of the type
$$[\pi(a_0') [D,\pi(a_1')][D,\pi(a_2')], \pi(a_3')].$$
Using the previous result, we can decompose it into the sum of
commutators of homogeneous polynomials and we immediately see that
for the algebra of the noncommutative torus commutators cannot
give a polynomial of degree $0$. Hence, no trivial cocycle can have
$\gamma$ as its image.

On the other hand, using the results on the Hochschild homology
of the noncommutative torus \cite{Wambst} we explicitly verify
that for the unique nontrivial cocycle (up to multiplication):
$$ c_0 = U^*V^* \ts V \ts U - V^* U^* \ts U \ts V, $$
its image, $\pi(c_0)$ vanishes. Therefore, in this case, the
Hochschild cocycle condition cannot be satisfied.

In the remaining case ($\phi=0$ and $\psi\not=0$, for instance) we use
similar arguments. We have:
$$ [D,\pi(U)]_\pm = \tau_\mu^\pm U I_r, \;\;\;
[D,\pi(V)]_\pm = (e^{\mp 2i \psi} -1) U D'_\pm, $$
where $I_r e_{\mu\nu,\pm} = e_{\mu\nu,\mp}$ and
$D' e_{\mu,nu,\pm} = e^{-2i\psi\nu} e_{\mu,nu,\mp}$.

Repeating the arguments from previous considerations, we obtain that
for any expression of the considered type, it might have three
components:
$$
\begin{aligned}
C(a_1,a_2) \pi(a_0)\pi(a_1)\pi(a_2) &+
C'(a_1,a_2) \pi(a_0)\pi(a_1)\pi(a_2) I_r D' \\
&+ C''(a_1,a_2) \pi(a_0)\pi(a_1)\pi(a_2) D' I_r,
\end{aligned}
$$
with three complex coefficients, depending only on the
multi-degree of $a_1,a_2$ and $\psi$.

The difference here is the appearance of $D'$, but again, it is sufficient
to verify that no such finite sum can be proportional to the identity, when
restricted to $\CV_+$, unless all $C',C''$ vanish and the degree of
$a_0 a_1 a_2$ is zero.

We can further follow the same arguments for the commutator presentation
of trivial Hochschild cycles, checking again explicitly that for the
nontrivial $c_0$, its image has non-zero coefficients $C'$ and $C''$ and
for this reason its image cannot be $\gamma$.

For the $\psi=0,\phi=0$ case, we calculate that
$$
\begin{aligned}
\gamma = \frac{1}{\tau_\mu^* \tau_\nu - \tau_\mu \tau_\nu^*}
& \left( \pi(V^*) \pi(U^*) [D,\pi(U)] [D,\pi(V)] \right. \\
& \left. -\pi(U^*) \pi(V^*) [D,\pi(V)] [D,\pi(U)] \right).
\end{aligned}
$$
provided that $\tau_\mu \tau_\nu^* \not= \tau_\mu^* \tau_\nu$.
\end{proof}

We can now state:
\begin{theorem}
There are four inequivalent equivariant spin structures on the
$2$-dimensional noncommutative torus, with a unique choice of
equivariant Dirac operator for each spin structure:
\begin{equation}
d_{\mu,\nu}^+ = \tau_\mu \mu + \tau_\nu \nu,
\label{dirac}
\end{equation}
which satisfies the Hochschild cycle condition, provided that
$\tau_\mu \tau_\nu^* \not= \tau_\mu^* \tau_\nu$. The spectrum
of the equivariant Dirac Operator depends on the spin structure.
\end{theorem}
\begin{proof}

As the previous lemmas showed the construction of the spectral data,
we only need to show their inequivalence and the dependence of the
spectrum of the Dirac operator on the choice of the class.

In order to see that the different reality structures we found
are not equivalent we need to find the same presentation of the
spectral geometries. This is achieved by relabelling the indices
so that their are all integers.

We obtain, on the Hilbert space with the basis labelled by
integers $m,n$:
\begin{center}
\begin{tabular}{|l|l||c|c|}
\hline
$\epsilon_\mu$ & $\epsilon_\nu$ & $J \ket{m,n,\pm}$ & $D \ket{m,n,\pm}$ \\ \hline
$0$ & $0$   &     $\pm \lambda^{-mn}\ket{-m,-n,\mp}$
& $(\tau_1 m + \tau_2 n) \ket{-m,-n,\pm}$ \\ \hline
$0$ & $\oh$ &   $\pm \lambda^{-m(n \mp \oh)} \ket{-m,-n+1,\mp}$
& $(\tau_1 m + \tau_2 n + \oh) \ket{-m,-n,\pm}$ \\ \hline
$\oh$ & $0$ &   $\pm \lambda^{-(m \mp \oh)n} \ket{-m+1,-n,\mp}$
& $(\tau_1 m + \tau_2 n +\oh) \ket{-m,-n,\pm}$ \\ \hline
$\oh$ & $\oh$ & $\pm \lambda^{-(m \mp \oh)(n \mp \oh)} \ket{-m+1,-n+1,\mp}$
& $(\tau_1 m + \tau_2 n + 1) \ket{-m,-n,\pm}$ \\ \hline
\end{tabular}
\end{center}
So two of the above cases have clearly a different spectrum
of the Dirac operator than the other two (for instance, for
$\epsilon_\mu \not=\epsilon_\nu$, $0$ is not in the spectrum
of $D$). It is thus immediately evident that the two pairs
corresponding to the different spectra of $D$ are not unitarily
equivalent.\\

Let us prove the mutual inequivalence of the two cases within
each pair. We shall show that there does not exist a unitary
operator $W$ on $\CV$, such that
$J_{\epsilon_\mu'\epsilon_\nu'} = WJ_{\epsilon_\mu\epsilon_\nu}W^*$,
and which leaves the remaining data of the spectral triple unchanged,
in particular
\begin{eqnarray}
 W\gamma W^* & = & \gamma \label{gamm} \\
 W\pi(a)W^* & = &  \pi(a) \qquad\qquad \forall a \in \CA(T_\lambda) \label{alg}
\end{eqnarray}

{}From the first of these equations, (\ref{gamm}), it follows that $W$ is blockdiagonal, $W:\CV_\pm \to \CV_\pm$.
To make use of (\ref{alg}) we first observe that
$e_{n,m,\pm} =  \pi(U^n) \pi(V^m) e_{0,0,\pm}$ for all $n,m$.
We shall denote
$$
W e_{0,0,\pm} = \sum_{k,l} w^\pm_{kl}\, e_{k,l,\pm}
$$
Using (\ref{alg}), we have
\begin{equation}
\begin{aligned}
W e_{n,m,\pm} & = &  W \pi(U^n) \pi(V^m)\, e_{0,0,\pm}
=  \pi(U^n) \pi(V^m) W\, e_{0,0,\pm} \\
& = & \pi(U^n) \pi(V^m) \left(\sum_{k,l} w^\pm_{kl}\, e_{k,l,\pm} \right) \\
& = & \sum_{k,l} \lambda^{-mk} w^\pm_{kl}\, e_{k+n,l+m,\pm}.
\end{aligned}
\end{equation}
Thus $W$ is completely determined by the coefficients $w^\pm_{kl}$.
The requirement $W J_{\epsilon_\mu\epsilon_\nu}W^* =J_{\epsilon_\mu'\epsilon_\nu'}$,
gives the following equation:
\begin{equation}
w^\pm_{-k,-l} \, j_{\epsilon_\mu, \epsilon_\nu,\pm}(n+k,m+l) =
\lambda^{2mk}\, w^\mp_{k,l} \, j_{\epsilon_\mu', \epsilon_\nu',\pm}(n,m),
\;\;\; \forall n,m,k,l \in \Z,
\end{equation}
where
$$ J_\iota\, e_{n,m,\pm}  = j_\iota(n,m)\, e_{-n,-m,\mp}, \;\;\;
  \iota= \{(\epsilon_\mu\epsilon_\nu),(\epsilon_\mu'\epsilon_\nu')\}.
$$
Taking the  $j_\iota(n,m)$ from the table above and inserting them into
the equation for the $w^\pm_{kl}$ one easily sees that there only exists
a solution if $J_{\epsilon_\mu',\epsilon_\nu'} =J_{\epsilon_\mu,\epsilon_\nu}$
in which case $W=\pm\id$. Thus there does not exist a unitary that intertwines
distinct reality structures $J_{\epsilon_\mu\epsilon_\nu}$.
\end{proof}

Finally, it is interesting to note that, were $W$ not required to
commute with $\gamma$ and the algebra representation, there would exist
such a unitary. For example, an unitary $W$ which intertwines
$J_{00}$ and $J_{0\oh}$,
$$ J_{0\oh} = W^* J_{00} W, $$
is given by:
$$ W e_{m,n,+} = \lambda^{\oh m} e_{m,n-1,-},
\;\;\;\; W e_{m,n,-} = e_{m,n,+}.$$
Indeed, then
\begin{equation}
\begin{split}
W^* J_{00} W e_{m,n,+} &= W^* J_{00} \lambda^{\oh m}\, e_{m,n-1,-} = \\
&= W^* \lambda^{-m(n-1)} \lambda^{-\oh m} \,e_{-m,-n+1,+} \\
&= \lambda^{-m(n-\oh)} \,e_{-m,-n+1,-},
\end{split}
\end{equation}
and
\begin{equation}
\begin{split}
W^* J_{00} W e_{m,n,-} &= W^* J_{00} e_{m,n,+}
= W^* \lambda^{-mn} e_{-m,-n,-} \\
&= \lambda^{-\oh m} \lambda^{-mn} e_{-m,-n+1,+} \\
&= \lambda^{-m(n+\oh)} e_{-m,-n+1,+}.
\end{split}
\end{equation}

\section{Conclusions}

We have shown that the noncommutative torus has, similarly as in the classical
$\lambda=1$ case, four inequivalent spin structures. It is not surprising that
the spin structures are closely related to the reality structure $J$. It is quite instructive, however,
that the pure algebraic conditions for
$J$ and the Dirac operator are not sufficient and one needs either the Hochschild
cocycle condition or the restriction due to the spectral properties of the Dirac
operator. Note that in \cite{MP} this was also shown to rule out
spin bundles with the wrong topology over the commutative sphere $S^2$, leaving
precisely one (real) $Spin$ structure in that case.

The "nonexisting" spurious classes of the reality operator $J$, which do not
lead to true Dirac operators have no classical (commutative) counterpart. For
this reason, it is hard to compare the construction with the steps of
Connes' reconstruction theorem for spin geometries \cite{GBVF} in order to
see whether their existence is a shadow of some other structures.

\end{document}